\documentclass[a4paper,11pt]{article}
\usepackage{amssymb,amsmath}

\begin{document}
\newcommand{\PMod}[1]{\!\!\pmod{#1}}

\newcommand{\Legendre}[2]{(\leavevmode\kern.1em \raise.5ex\hbox{\the\scriptfont0 #1}\kern-.15em
/\kern-.15em\lower.45ex\hbox{\the\scriptfont0 #2})}
\newcommand{\LLegendre}[2]{\big(\leavevmode\kern.05em \raise.5ex\hbox{#1}\kern-.2em
/\kern-.15em\lower.55ex\hbox{#2}\big)}

\parindent = 0 pt

\newtheorem{theorem}{Theorem}
\newtheorem{conjecture}{Conjecture}

\begin{center}
{ \bf
 \Large{A LLT-like test for proving the primality of Fermat
 numbers.}
}

\vspace{.1in}

Tony Reix (Tony.Reix@laposte.net)

\vspace{.08in}

First version: 2004, 24th of September

Updated: 2005, 19th of October

\end{center}

\noindent

\vspace{.08in}

In 1876, \'{E}douard Lucas discovered a method for proving that a
number is prime or composite without searching its factors. His
method was based on the properties of the \emph{Lucas Sequences}. He
first used his method for Mersenne numbers and proved that
$2^{127}-1$ is a prime. In 1930, Derrick Lehmer provided a complete
and clean proof. This test of primality for Mersenne numbers is now
known as: Lucas-Lehmer Test (LLT).

\vspace{.1in}

Few people know that Lucas also used his method for proving that a
Fermat number is prime or composite, still with an unclear proof. He
used his method for proving that $2^{2^6}+1$ is composite. Lehmer
did not provide a proof of Lucas' method for Fermat numbers.

\vspace{.1in}

This paper provides a proof of a LLT-like test for Fermat numbers,
based on the properties of Lucas Sequences and based on the method
of Lehmer. The seed (the starting value $S_{\,0}$ of the
$\{S_{\,i}\}$ sequence) used here is 5, though Lucas used 6.

\vspace{.1in}

Primality tests for special numbers are classified into $N-1$ and
$N+1$ categories, meaning that the numbers $N-1$ or $N+1$ can be
completely or partially factored. Since many books talk about the
LLT only in the $N+1$ chapter for Mersenne numbers $N=2^q-1$, it
seemed useful to remind that the LLT can also be used for numbers
$N$ such that $N-1$ is easy to factor, like Fermat numbers
$N=2^{2^n}+1$, by providing a proof \emph{\`{a}la} Lehmer.

\begin{theorem}
\ \newline $F_n = 2^{\scriptstyle 2^{\scriptstyle n}} + 1$ ($n
\geqslant 1$) is a prime if and only if it divides $S_{\,
2^{\scriptstyle n}-2}$ , where \ $S_{\,0} = 5$ \ and \ $S_{\,i} =
S_{\,i-1}^{\,\,2} - 2$ \ for \ $i=1,2,3, ... \ 2^{\scriptstyle n}-2$
.
\end{theorem}
The proof is based on chapters 4 (The Lucas Functions) and 8.4 (The
Lehmer Functions) of the book "\'{E}douard Lucas and Primality
Testing" of H. C. Williams (A Wiley-Interscience publication, 1998).

\vspace{.1in}

Chapter 1 explains how the $(P,Q)$ parameters have been found. Then
Chapter 2 provides the Lehmer theorems used for the proof. Then
Chapter 3 and 4 provide the proof for: $F_n$ prime $\Longrightarrow
F_n \mid S_{2^{\scriptstyle n}-2}$ and the converse, proving theorem
1. Chapter 5 provides numerical examples. The appendix in Chapter 6
provides first values of $U_n$ and $V_n$ plus some properties.

\vspace{.15in}

century).

\newpage

\section{ Lucas Sequence with $P=\sqrt{R}$ }

Let $S_0=5$ and $S_i = S_{i-1}^2-2$ . $S_1=23$,
$S_2=527=17\times31$, ...

\vspace{-.2in}
$$
\text{It has been checked that: }\ \ \left\{
\begin{array}{lllcr}
S_{2^{\scriptstyle n}-2} \equiv 0 \ \PMod{F_n} \text{ \ for } n=1 ... 4 \\
S_{2^{\scriptstyle n}-2} \neq   0 \ \PMod{F_n} \text{ \ for } n=5 ... 14 \\
\end{array}
\right.
$$

Here after, we search a Lucas Sequence $(U_m)_{\,m \geqslant 0}$ and
its companion $(V_m)_{\,m\geqslant0}$ with $(P,Q)$ that fit with the
values of the $S_i$ sequence.

\vspace{.1in}

We define the Lucas Sequence $V_m$ such that:
\vspace{-.1in}
\begin{equation} \label{A}
V_{2^{\scriptstyle k+1}} = S_k
\end{equation}
\vspace{-.1in}
$$
\text{Thus we have: }\ \ \left\{
\begin{array}{lllcr}
V_2 & = & S_0 & = & 5 \\
V_4 & = & S_1 & = & 23 \\
V_8 & = & S_2 & = & 527
\end{array}
\right.
$$
$$
\text{If (4.2.7) page 74 ( $V_{2n} = V_n^2 - 2Q^n$ ) applies, we
have:} \ \ \left\{
\begin{array}{llrl}
V_4 = V_2^2 -2Q^2 \\
\vspace{-.15in} \\
V_8 = V_4^2 -2Q^4
\end{array}
\right.
$$
and thus: $Q = \sqrt[2]{\frac{V_2^2-V_4}{2}}=
\sqrt[4]{\frac{V_4^2-V_8}{2}}= \pm 1$ .

\vspace{.15in}

With (4.1.3) page 70 ( $V_{n+1} = PV_{n}-QV_{n-1}$ ), and with:
$$
\left\{
\begin{array}{llrl}
V_0 = 2 \\
V_1 = P \\
V_2 = PV_1 -QV_0 = P^2 -2Q
\end{array}
\right.
$$

we have: \ $P=\sqrt{V_2 + 2Q} = \sqrt{7}$ \ or \ $\sqrt{3}$ .

\vspace{.05in}

In the following we consider: $(P,Q) = (\sqrt{7},1)$ .

\vspace{.15in}

As explained by Williams page 196, "all of the identity relations
[Lucas functions] given in (4.2) continue to hold, as these are true
quite without regard as to whether $P, Q$ are integers".

\vspace{.1in}

So, like Lehmer, we define $P=\sqrt{R}$ such that $R=7$ and $Q=1$
are coprime integers and we define (Property (8.4.1) page 196):

\vspace{-.18in}
$$
\overline{V}_n = \ \left\{
\begin{array}{ll}
V_n & \text{ when } 2 \mid n \\
\vspace{-.1in} \\
V_n / \sqrt{R} & \text{ when } 2 \nmid n
\end{array}
\right.
\ \ \ \ \ \overline{U}_n = \ \left\{
\begin{array}{ll}
U_n / \sqrt{R} & \text{ when } 2 \mid n \\
\vspace{-.1in} \\
U_n & \text{ when } 2 \ \nmid n
\end{array}
\right.
$$
\vspace{-.02in} in such a way that $\overline{V}_n$ and
$\overline{U}_n$ are always integers.

\vspace{.1in}

Tables 1 to 5 give values of $U_i$ , $V_i$ , $\overline{U}_i
\PMod{F_n}$ , $\overline{V}_i \PMod{F_n}$ , with $(P,Q) =
(\sqrt{7},1)$ , for $n = 1, 2, 3, 4$ .

\section{ Lehmer theorems }
Like Lehmer, let define the symbols (where \LLegendre{a\,}{\,b} is
the Legendre symbol): \vspace{-.1in}
$$
\left\{
\begin{array}{l}
\varepsilon = \varepsilon(p) = \LLegendre{D}{p} \\
\vspace{-.1in} \\
\sigma = \sigma(p) = \LLegendre{R}{p} \\
\vspace{-.1in} \\
\tau = \tau(p) = \LLegendre{Q}{p}
\end{array}
\right.
$$
\vspace{-.1in}

The 2 following formulas (from page 77) will help proving
properties: \vspace{-.1in}
$$
\text{(4.2.28) \ \ \ } 2^{m-1}U_{mn} \ = \ \sum_{i=0}^{\lfloor
\frac{\scriptstyle m}{\scriptstyle 2} \rfloor} {m \choose 2i+1} D^i
U^{2i+1}_n V^{m-(2i+1)}_n \
$$
\vspace{-.1in}
$$
\text{(4.2.29) \ \ \ } 2^{m-1}V_{mn} \ = \ \ \sum_{i=0}^{\lfloor
\frac{\scriptstyle m}{\scriptstyle 2} \rfloor} {m \choose 2i} D^i
U^{2i}_n V^{m-2i}_n \ \ \ \ \ \ \ \ \ \ \ \ \ \ \ \
$$

\textbf{Property (8.4.2)} page 196 : \vspace{-.2in}
$$
\text{If } p \text{ is an odd prime and } p \nmid Q \text{, then: }
\left\{
\begin{array}{l}
\overline{U}_{p} \equiv \LLegendre{D}{p} \PMod{p} \\
\vspace{-.1in} \\
\overline{V}_{p} \equiv \LLegendre{R}{p} \PMod{p}
\end{array}
\right.
$$

\vspace{-.06in}

\textbf{Proof:}

Since $p$ is a prime, and by Fermat little theorem, we have:
$2^{p-1} \equiv 1 \ \PMod{p}$.

\vspace{.05in}

$\bullet $ By (4.2.28), with $m=p$ and $n=1$, since $U_1 = 1$ and
$V_1=P$, we have: \vspace{-.1in}
$$
2^{p-1}U_{p} = \ \sum_{i=0}^{ \frac{\scriptstyle p-1}{\scriptstyle
2}} {p \choose 2i+1} D^i U^{2i+1}_1 V^{p-(2i+1)}_1
$$
$$
2^{p-1}U_{p} = {p \choose 1} P^{p-1} + {p \choose
3} DP^{p-3} + ... + {p \choose p} D^{\frac{\scriptstyle
p-1}{\scriptstyle 2}}P^0
$$
Since ${p \choose i} \equiv 0 \ \PMod{p}$ when $0 < i < p$ and ${p
\choose p} = 1$ , we have: \vspace{-.1in}
$$
U_{p} = \overline{U}_{p} \equiv D^{\frac{\scriptstyle
p-1}{\scriptstyle 2}} \equiv \LLegendre{D}{p} \ \PMod{p}
$$

$\bullet $ By (4.2.29), with $m=p$ and $n=1$, since $U_1 = 1$ and
$V_1=P$, we have:  \vspace{-.1in}
$$
2^{p-1}V_{p} = \ \sum_{i=0}^{\frac{\scriptstyle p-1}{\scriptstyle
2}} {p \choose 2i} D^i U^{2i}_1 V^{p-2i}_1
$$
$$
2^{p-1}V_{p} = {p \choose 0} P^{p} + {p \choose 2} DP^{p-2} + ... +
{p \choose p-1} D^{\frac{\scriptstyle p-1}{\scriptstyle 2}}P
$$
Since ${p \choose 0} = 1$ , and ${p \choose i} \equiv 0 \ \PMod{p}$
when $0 < i < p$ , we have:  \vspace{-.1in}
$$
V_{p} \equiv P^{p} \ \text{ and } \ \overline{V}_{p} \equiv P^{p-1}
\equiv R^{\frac{\scriptstyle p-1}{\scriptstyle 2}} \equiv
\LLegendre{R}{p} \ \PMod{p}
$$
\vspace{-.45in}
\begin{flushright}
$\square$
\end{flushright}

\textbf{Property (8.4.3)} page 197 : \vspace{-.05in}
\begin{center}
$p$ odd prime and $p \nmid Q$ $\Longrightarrow p \mid
\overline{U}_{\scriptstyle p-\sigma\varepsilon}$
\end{center}

\vspace{.05in}

\textbf{Proof}

By (4.2.28) with $n=1$, $V_1=P$, since $p$ is a prime and $(R,Q) =
1$, we have:

\vspace{.1in}

$\bullet \ \ \text{With: } m=p+1$ \vspace{-.15in}
$$
2^pU_{p+1} = \ \sum_{i=0}^{ \frac{\scriptstyle p+1}{\scriptstyle 2}}
{p+1 \choose 2i+1} D^i P^{p-2i}
$$
$$
2^pU_{p+1} = {p+1 \choose 1} P^p + {p+1 \choose 3} D P ^{p-2} + ...
+ {p+1 \choose p}D^{\frac{\scriptstyle p-1}{\scriptstyle 2}} P +
{p+1 \choose p+2}D^{\frac{\scriptstyle p+1}{\scriptstyle 2}} P^{-1}
$$
$$
2^pU_{p+1} = (p+1)P^p + (p+1)p\big[...\big] + (p+1)D^{\frac{p-1}{2}}
P + 0 D^{\frac{p+1}{2}} P^{-1}
$$
$$
2^pU_{p+1} = P^p + D^{\frac{p-1}{2}}P + p\big[...\big] =
P\big[(P^2)^{\frac{p-1}{2}}+D^{\frac{p-1}{2}}\big] + p\big[...\big]
$$
$$
\frac{2^pU_{p+1}}{P} = 2^p\overline{U}_{p+1} \equiv
R^{\frac{p-1}{2}}+D^{\frac{p-1}{2}} \equiv \LLegendre{R}{p} +
\LLegendre{D}{p} = \sigma(p) + \varepsilon(p) \PMod{p}
$$

\begin{center}
Thus, if $\sigma \varepsilon = \sigma(p) \times \varepsilon(p) = -1$
, then $p \mid \overline{U}_{p+1} = \overline{U}_{p-\sigma\epsilon}$
.
\end{center}

$\bullet \ \ \text{With: } m=p-1$ :  \vspace{-.15in}
$$
2^{p-2}U_{p-1} = \ \sum_{i=0}^{ \frac{\scriptstyle p-1}{\scriptstyle
2}} {p-1 \choose 2i+1} D^i P^{p-2(i+1)}
$$
$$
2^{p-2}U_{p-1} = {p-1 \choose 1} P^{p-2} + {p-1 \choose 3} D P
^{p-4} + ... + {p-1 \choose p-2}D^{\frac{\scriptstyle
p-3}{\scriptstyle 2}} P + {p-1 \choose p}D^{\frac{\scriptstyle
p-1}{\scriptstyle 2}} P^{-1}
$$
$$
2^{p-2}U_{p-1} = (p-1)P^{p-2} + (p-1)DP^{p-4} + ... + (p-1)
D^{\frac{\scriptstyle p-3}{\scriptstyle 2}} P + 0
D^{\frac{\scriptstyle p-1}{\scriptstyle 2}} P^{-1}
$$
$$
\frac{2^{p-2}U_{p-1}}{P} \equiv -[P^{p-3} + DP^{p-5} + ... \ +
D^{\frac{\scriptstyle p-3}{\scriptstyle 2}}] \equiv - \frac{P^{p-1}
- D^{\frac{\scriptstyle p-1}{\scriptstyle 2}}}{P^2-D} \PMod{p}
$$
$$
2^{p-2}\overline{U}_{p-1}(P^2-D) \equiv - (P^2)^{\frac{\scriptstyle
p-1}{\scriptstyle 2}} + D^{\frac{\scriptstyle p-1}{\scriptstyle 2}}
\equiv \varepsilon(p) - \sigma(p) \PMod{p}
$$

\begin{center}
Thus, if $ \sigma \varepsilon = \sigma(p) \times \varepsilon(p) = 1$
, then $p \mid \overline{U}_{p-1} = \overline{U}_{p-\sigma\epsilon}$
.
\end{center}
\vspace{-.3in}
\begin{flushright}
$\square$
\end{flushright}

\vspace{-.05in}

\textbf{Property (8.4.4) page 197} \vspace{-.15in}
$$
\text{If } p \text{ is an odd prime and } p \nmid Q \text{, then: }
V_{p-\sigma \varepsilon} \equiv 2\sigma Q^{\frac{\scriptstyle
1-\sigma \varepsilon}{\scriptstyle 2}} \PMod{p} \ .
$$

\vspace{.05in}

\begin{theorem}[8.4.1]
If $p$ is an odd prime and $p \nmid QRD$ , then: \vspace{-.05in}
$$
\left\{
\begin{array}{lll}
p \mid \overline{V}_{\frac{\scriptstyle p -
\sigma\epsilon}{\scriptstyle 2}} & \text{ when} & \sigma = -\tau
\\
\vspace{-.1in} \\
p \mid \overline{U}_{\frac{\scriptstyle p -
\sigma\epsilon}{\scriptstyle 2}} & \text{ when} & \sigma = \tau
\\
\end{array}
\right.
$$
\end{theorem}

\vspace{.05in}

\textbf{Definition (8.4.2) page 197 of $\omega(m)$ : } For a given
$m$, denote by $\omega = \omega(m)$ the value of the least positive
integer $k$ such that $m \mid \overline{U}_k$ . If $\omega(m)$
exists, $\omega(m)$ is called the \textbf{rank of apparition} of $m$
.

\vspace{.05in}

\begin{theorem}[8.4.3] $$
\left\{
\begin{array}{l}
\text{If } k \mid n \text{, then } \overline{U}_k \mid
\overline{U}_n \ . \\
\vspace{-.1in} \\
\text{If } m \mid \overline{U}_n \text{, then } \omega(m) \mid n \ .
\end{array}
\right.
$$
\end{theorem}

\vspace{.02in}

\begin{theorem}[8.4.5]
\ \ If $(m,Q) = 1$ , then $\omega(m)$ exists.
\end{theorem}

\vspace{.02in}

\begin{theorem}[8.4.6]
\ \ \ If $(N,2QRD) = 1$ and $N \pm 1$ is the rank of apparition of
$N$, then $N$ is a prime.
\end{theorem}

\vspace{.02in}

\begin{theorem}[8.4.7]
\ \ \ If $(N,2QRD) = 1$ , $\overline{U}_{N \pm 1} \equiv 0 \
\PMod{N}$ and $\overline{U}_{\scriptstyle \frac{N \pm
1}{\scriptstyle q}} \neq 0 \ \PMod{N}$ for each distinct prime
divisor $q$ of $N \pm 1$, then $N$ is a prime.
\end{theorem}

\textbf{Proof:}

Let $\omega = \omega(N)$ . We see that $\omega \mid N \pm 1$ , but
$\omega \nmid (N \pm 1)/q$ . Thus if $q^\alpha \parallel N \pm 1$ ,
then $q^\alpha \mid \omega$ . It follows that $\omega = N \pm 1$ and
$N$ is a prime by Theorem 5 (8.4.6) .

\vspace{.1in}

\section{ $F_n$ prime $\Longrightarrow  F_n \mid
\overline{V}_\frac{\scriptstyle F_{\scriptstyle n}-1}{2}$ and $F_n
\mid S_{2^{\scriptstyle n}-2}$ }

Let $N = F_n = 2^{2^{\scriptstyle n}}+1$ with $n \geq 1$ be an odd
prime.

\vspace{.05in}

Let: $P = \sqrt{R}$ , $R = 7$ , $Q = 1$ , and $D = P^2 - 4Q = 3$ .

\vspace{.12in}

Hereafter we compute $\LLegendre{3}{N}$ and $\LLegendre{7}{N}$ :

\vspace{.1in}

$\bullet \ \LLegendre{3}{N} : $ \vspace{-.35in}
$$
\text{    Since:  } \left\{
\begin{array}{l}
N \ \ \text{odd prime} \\
\vspace{-.18in} \\
N = (4)^{2^{\scriptstyle n-1}}+1 \equiv 2 \PMod{3} \\
\vspace{-.15in} \\
\LLegendre{N}{3} = \LLegendre{2}{3} = -1 \\
\vspace{-.15in} \\
 \LLegendre{3}{N} = \LLegendre{N}{3} \times
(-1)^{\frac{\scriptstyle 3-1}{\scriptstyle 2}\frac{\scriptstyle
N-1}{\scriptstyle 2}}
\end{array}
\right. \text{ \ \ \  then:   } \LLegendre{3}{N} = -1 \ .
$$

$\bullet \ \LLegendre{7}{N} :$ We have: \vspace{-.20in}
$$
\left\{
\begin{array}{llll}
2^3      & \equiv & 1   & \PMod{7} \\
\vspace{-.15in} \\
2^{3a+b} & \equiv & 2^b & \PMod{7}
\end{array}
\right.
$$

\vspace{-.05in}

With $2^n \equiv b \ \PMod{3}$ , we have: $2^{2^{\scriptstyle n}}+1
\equiv 2^{b} + 1 \ \PMod{7}$ . Then we study the exponents of 2,
modulo 3 . We have: $2^2 \equiv 1 \ \PMod{3}$ , and: \vspace{-.05in}
$$
\text{If } n = 2m \ \ \ \ \ \, \left\{
\begin{array}{llll}
2^{2m} \equiv 1 \PMod{3} \\
\vspace{-.15in} \\
N=2^{2^{\scriptstyle 2m}}+1 \equiv 2^1 + 1 \equiv 3 \PMod{7} \\
\vspace{-.1in} \\
\LLegendre{N}{7} = \LLegendre{3}{7} = -1
\end{array}
\right.
$$

$$
\text{ \ \ If } n = 2m+1 \ \left\{
\begin{array}{llll}
2^{2m+1} \equiv 2 \PMod{3} \\
\vspace{-.15in} \\
N=2^{2^{\scriptstyle 2m+1}}+1 \equiv 2^2 + 1 \equiv 5 \PMod{7} \\
\vspace{-.1in} \\
\LLegendre{N}{7} = \LLegendre{5}{7} = -1
\end{array}
\right.
$$

\ \ \ \ \ \ Finally, we have: $\LLegendre{7}{N} = \LLegendre{N}{7}
(-1)^{\frac{\scriptstyle 7-1}{\scriptstyle 2} 2^{\scriptstyle
2^{\scriptstyle n}}} = \LLegendre{N}{7} = -1$ .

\vspace{.1in}
$$
\text{So we have: \ }\left\{
\begin{array}{lr}
\varepsilon = \LLegendre{D}{N} = \LLegendre{3}{N} = & -1 \\
\vspace{-.1in} \\
\sigma = \LLegendre{R}{N} = \LLegendre{7}{N} = & -1 \\
\vspace{-.1in} \\
\tau = \LLegendre{Q}{N} = \LLegendre{1}{N} = & + 1
\end{array}
\right.
$$
\vspace{.05in}

Since $\sigma = - \tau$ , $\sigma \epsilon = +1$ , and $F_n \nmid
QRD$ with $n \geq 1$, then by Theorem 2 (8.4.1) we have:
$$F_n \text{ prime }  \Longrightarrow F_n \mid \overline{V}_{\frac{\scriptstyle F_n-1}{\scriptstyle 2}} =
V_{2^{\scriptstyle 2^{\scriptstyle n}-1}}$$

By (1) we have: $V_{\scriptstyle 2^{{\scriptstyle k}-1}} = S_{k-2}$
and thus, with $k=2^n$: $F_n \mid S_{2^{\scriptstyle n}-2}$ .
\vspace{-.1in}
\begin{flushright}
$\square$
\end{flushright}

\vspace{-.2in}

\section{ $F_n \mid S_{2^{\scriptstyle n}-2} \ \Longrightarrow \ F_n$ is a prime }

Let $N = F_n$ with $n \geq 1$ . By (1) we have: $N \mid
S_{2^{\scriptstyle n}-2} \Longrightarrow N \mid V_{2^{\scriptstyle
2^{\scriptstyle n}-1}}$ .

\vspace{.03in}

And thus, by (4.2.6) page 74 \ ( $U_{2a} = U_a V_a$ ) , we have: $N
\mid \overline{U}_{2^{\scriptstyle 2^{\scriptstyle n}}}$ .

\vspace{.15in}

By (4.3.6) page 85: \ ( $(V_n,U_n) \mid 2Q^{\,n}$ for any $n$ ), and
since $Q = 1$ , then:

\vspace{.04in}

$(V_{2^{\scriptstyle 2^{\scriptstyle n}-1}} ,
\overline{U}_{2^{\scriptstyle 2^{\scriptstyle n}-1}}) = 2$ and thus:
$N \nmid \overline{U}_{2^{\scriptstyle 2^{\scriptstyle n}-1}}$ since
$N$ odd.

\vspace{.15in}

With $\omega = \omega(N)$ , by Theorem 3 (8.4.3) we have : $\omega
\mid 2^{2^{\scriptstyle n}}$ and $\omega \nmid 2^{2^{\scriptstyle
n}-1}$ .

\vspace{.05in}

This implies: $\omega = 2^{2^{\scriptstyle n}} = N-1$ . Then $N-1$
is the rank of apparition of N, and thus by Theorem 5 (8.4.6) N is a
prime. \vspace{-.1in}
\begin{flushright}
$\square$
\end{flushright}

\vspace{.05in}

This test of primality for Fermat numbers has been communicated to
the community of number theorists working on this area on
mersenneforum.org
(http://www.mersenneforum.org/showthread.php?t=2130) in May 2004,
and the proof was finalized in September 2004.

Then, in a private communication, Robert Gerbicz provided a proof of
the same theorem based on $Q[\sqrt{21}]$.

\section{ Numerical Examples }

$\PMod{F_2}$ $S_0 = 5 \stackrel{1}{\mapsto} {6}
\stackrel{2}{\mapsto} S_{2^2-2} \equiv 0 \ $

$\PMod{F_3}$ $S_0 = 5 \stackrel{1}{\mapsto} 23 \stackrel{2}{\mapsto}
13 \stackrel{3}{\mapsto} 167 \stackrel{4}{\mapsto} {131}
\stackrel{5}{\mapsto} {197} = -60 \stackrel{6}{\mapsto} S_{2^3-2}
\equiv 0$

$\PMod{F_4}$ $S_0 = 5 \stackrel{1}{\mapsto} 23 \stackrel{2}{\mapsto}
527 \stackrel{3}{\mapsto} {15579} \stackrel{4}{\mapsto} {21728}
\stackrel{5}{\mapsto} 42971 \stackrel{6}{\mapsto} 1864
\stackrel{7}{\mapsto} 1033 \stackrel{8}{\mapsto} 18495
\stackrel{9}{\mapsto} 27420 \stackrel{10}{\mapsto}15934
\stackrel{11}{\mapsto} 2016 \stackrel{12}{\mapsto} 960
\stackrel{13}{\mapsto} 4080 \stackrel{14}{\mapsto} S_{2^4-2} \equiv
0$

\section{ Appendix: Table of $U_i$ and $V_i$ }

With $n = 2, 3, 4$, we have the following (not proven) properties
(modulo $F_n$):
$$
\left\{
\begin{array}{llr}
\overline{U}_{F_{\scriptstyle n}-5} & \equiv & 5 \\
\vspace{-.15in} \\
\overline{U}_{F_{\scriptstyle n}-4} & \equiv & 6 \\
\vspace{-.15in} \\
\overline{U}_{F_{\scriptstyle n}-3} & \equiv & 1 \\
\vspace{-.15in} \\
\overline{U}_{F_{\scriptstyle n}-2} & \equiv & 1 \\
\vspace{-.15in} \\
\overline{U}_{F_{\scriptstyle n}-1} & \equiv & 0 \\
\vspace{-.15in} \\
\overline{U}_{F_{\scriptstyle n}}   & \equiv & -1 \\
\vspace{-.15in} \\
\overline{U}_{F_{\scriptstyle n}+1} & \equiv & -1 \\
\vspace{-.15in} \\
\overline{U}_{F_{\scriptstyle n}+2} & \equiv & -6 \\
\vspace{-.15in} \\
\overline{U}_{F_{\scriptstyle n}+3} & \equiv & -5
\end{array}
\right. \ \ \ \ \ \ \ \ \ \ \ \ \left\{
\begin{array}{llr}
\overline{V}_{F_{\scriptstyle n}-5} & \equiv & -23 \\
\vspace{-.15in} \\
\overline{V}_{F_{\scriptstyle n}-4} & \equiv & -4 \\
\vspace{-.15in} \\
\overline{V}_{F_{\scriptstyle n}-3} & \equiv & -5 \\
\vspace{-.15in} \\
\overline{V}_{F_{\scriptstyle n}-2} & \equiv & -1  \\
\vspace{-.15in} \\
\overline{V}_{F_{\scriptstyle n}-1} & \equiv & -2 \\
\vspace{-.15in} \\
\overline{V}_{F_{\scriptstyle n}}   & \equiv & -1 \\
\vspace{-.15in} \\
\overline{V}_{F_{\scriptstyle n}+1} & \equiv & -5 \\
\vspace{-.15in} \\
\overline{V}_{F_{\scriptstyle n}+2} & \equiv & -4 \\
\vspace{-.15in} \\
\overline{V}_{F_{\scriptstyle n}+3} & \equiv & -23 \\
\end{array}
\right.
$$

\begin{table}[htbp]
\begin{center}
\begin{tabular}{|r|rl|rl|}
\hline
$i$ & $U_i$ & & $V_i$ & \\
\hline
0  &             0 & $\times \sqrt{7}$ &              2 &                   \\
1  &             1 &                   &              1 & $\times \sqrt{7}$ \\
2  &             1 & $\times \sqrt{7}$ &              5 &                   \\
3  &             6 &                   &              4 & $\times \sqrt{7}$ \\
4  &             5 & $\times \sqrt{7}$ &             23 &                   \\
5  &            29 &                   &             19 & $\times \sqrt{7}$ \\
6  &            24 & $\times \sqrt{7}$ &            110 &                   \\
7  &           139 &                   &             91 & $\times \sqrt{7}$ \\
8  &           115 & $\times \sqrt{7}$ &            527 &                   \\
9  &           666 &                   &            436  & $\times \sqrt{7}$ \\
10 &           551 & $\times \sqrt{7}$ &           2525 &                   \\
11 &          3191 &                   &           2089  & $\times \sqrt{7}$ \\
12 &          2640 & $\times \sqrt{7}$ &          12098 &                   \\
13 &         15289 &                   &          10009  & $\times \sqrt{7}$ \\
14 &         12649 & $\times \sqrt{7}$ &          57965 &                   \\
15 &         73254 &                   &          47956  & $\times \sqrt{7}$ \\
16 &         60605 & $\times \sqrt{7}$ &         277727 &                   \\
17 &        350981 &                   &         229771  & $\times \sqrt{7}$ \\
18 &        290376 & $\times \sqrt{7}$ &        1330670 &                   \\
19 &       1681651 &                   &        1100899  & $\times \sqrt{7}$ \\
20 &       1391275 & $\times \sqrt{7}$ &        6375623 &                   \\
21 &       8057274 &                   &        5274724  & $\times \sqrt{7}$ \\
22 &       6665999 & $\times \sqrt{7}$ &       30547445 &                   \\
23 &      38604719 &                   &       25272721  & $\times \sqrt{7}$ \\
24 &      31938720 & $\times \sqrt{7}$ &      146361602 &                   \\
25 &     184966321 &                   &      121088881  & $\times \sqrt{7}$ \\
26 &     153027601 & $\times \sqrt{7}$ &      701260565 &                   \\
27 &     886226886 &                   &      580171684  & $\times \sqrt{7}$ \\
28 &     733199285 & $\times \sqrt{7}$ &     3359941223 &                   \\
29 &    4246168109 &                   &     2779769539  & $\times \sqrt{7}$ \\
30 &    3512968824 & $\times \sqrt{7}$ &    16098445550 &                   \\
31 &   20344613659 &                   &    13318676011  & $\times \sqrt{7}$ \\
32 &   16831644835 & $\times \sqrt{7}$ &    77132286527 &                   \\
33 &   97476900186 &                   &    63813610516  & $\times \sqrt{7}$ \\
34 &   80645255351 & $\times \sqrt{7}$ &   369562987085 &                   \\
35 &  467039887271 &                   &   305749376569  & $\times \sqrt{7}$ \\
36 &  386394631920 & $\times \sqrt{7}$ &  1770682648898 &                   \\
37 & 2237722536169 &                   &  1464933272329  & $\times \sqrt{7}$ \\
38 & 1851327904249 & $\times \sqrt{7}$ &  8483850257405 &                   \\
39 &10721572793574 &                   &  7018916985076  & $\times \sqrt{7}$ \\
40 & 8870244889325 & $\times \sqrt{7}$ & 40648568638127 &                   \\
\hline
\end{tabular}
\end{center}
\caption{  $P=\sqrt{7} \ \ , \ \ Q = 1$ } \label{tab:table1}
\end{table}

\begin{table}[htbp]
\begin{center}
\begin{tabular}{|r|r|r|}
\hline
$i$ & $\overline{U}_i \ \PMod{F_1}$ & $\overline{V}_i \ \PMod{F_1}$ \\
\hline
0  &         0  &   2 \\
1  &         1  &   1 \\
2  &         1  &   \textbf{0} \\
3  &         1  &   4 \\
4  & \textbf{0} &   3 \\
\textbf{5} & 4  &   4 \\
6  &         4  &   0 \\
7  &         4  &   1 \\
8  &         0  &   2 \\
\hline
\end{tabular}
\end{center}
\caption{  $P=\sqrt{7} \ \ , \ \ Q = 1$ , Modulo $F_1$ }
\label{tab:table2}
\end{table}

\begin{table}[htbp]
\begin{center}
\begin{tabular}{|r|r|r|}
\hline
$i$ & $\overline{U}_i \ \PMod{F_2}$ & $\overline{V}_i \ \PMod{F_2}$ \\
\hline
0  &         0  &   2 \\
1  &         1  &   1 \\
2  &         1  &   5 \\
3  &         6  &   4 \\
4  &         5  &   6 \\
5  &        12  &   2 \\
6  &         7  &   8 \\
7  &         3  &   6 \\
8  &          13 & \textbf{0} \\
9  &           3 & 11 \\
10 &           7 &  9 \\
11 &          12 & 15 \\
12 &           5 & 11 \\
13 &           6 & -4 \\
14 &           1 & -5 \\
15 &           1 & -1 \\
16 &  \textbf{0} & -2 \\
\textbf{17} & -1 & -1 \\
18 &          -1 & -5 \\
19 &          -6 & -4 \\
20 &          -5 & 11 \\
21 &           5 & 15 \\
22 &          10 &  9 \\
23 &          14 & 11 \\
24 &           4 &  0 \\
\hline
\end{tabular}
\end{center}
\caption{  $P=\sqrt{7} \ \ , \ \ Q = 1$ , Modulo $F_2$ }
\label{tab:table3}
\end{table}

\begin{table}[htbp]
\begin{center}
\begin{tabular}{|r|r|r|}
\hline
$i$ & $\overline{U}_i \ \PMod{F_3}$ & $\overline{V}_i \ \PMod{F_3}$ \\
\hline
0   &           0 &   2 \\
1   &           1 &   1 \\
2   &           1 &   5 \\
3   &           6 &   4 \\
4   &           5 &  23 \\
8   &         115 &  13 \\
16  &         210 & 167 \\
32  &         118 & 131 \\
64  &          38 & 197 \\
128 &          33 & \textbf{0} \\
192 &          38 &  60 \\
224 &         118 & 126 \\
240 &         210 &  90 \\
248 &         115 & -13 \\
252 &           5 & -23 \\
253 &           6 &  -4 \\
254 &           1 &  -5 \\
255 &           1 &  -1 \\
256 &  \textbf{0} &  -2 \\
\textbf{257} & -1 &  -1 \\
258 &          -1 &  -5 \\
259 &          -6 &  -4 \\
260 &          -5 & -23 \\
\hline
\end{tabular}
\end{center}
\caption{  $P=\sqrt{7} \ \ , \ \ Q = 1$ , Modulo $F_3$ }
\label{tab:table4}
\end{table}

\begin{table}[htbp]
\begin{center}
\begin{tabular}{|r|r|r|}
\hline
$i$ & $\overline{U}_i \ \PMod{F_4}$ & $\overline{V}_i \ \PMod{F_4}$ \\
\hline
 2048 &        9933 & 15934 \\
 4096 &         567 &  2016 \\
 8192 &       28943 &   960 \\
16384 &       63129 &  4080 \\
32768 &        5910 & \textbf{0} \\
65532 &           5 &  -23 \\
65533 &           6 &   -4 \\
65534 &           1 &   -5 \\
65535 &           1 &   -1 \\
65536 &  \textbf{0} &   -2 \\
\textbf{65537} & -1 &   -1 \\
65538 &          -1 &   -5 \\
65539 &          -6 &   -4 \\
65540 &          -5 &  -23 \\
\hline
\end{tabular}
\end{center}
\caption{  $P=\sqrt{7} \ \ , \ \ Q = 1$ , Modulo $F_4$ }
\label{tab:table5}
\end{table}

The values of $\overline{U'}_n$ and $\overline{V'}_n$ $(_{n \geq
1})$ with $(P,Q) = (\sqrt{3},-1)$ can be built by:
$$
\left\{
\begin{array}{lll}
\overline{U'}_{2n} & = & \overline{U}_{2n} \\
\vspace{-.1in} \\
\overline{U'}_{2n+1} & = & \overline{V}_{2n+1}
\end{array}
\right. \ \ \ \left\{
\begin{array}{lll}
\overline{V'}_{2n} & = & \overline{V}_{2n} \\
\vspace{-.1in} \\
\overline{V'}_{2n+1} & = & \overline{U}_{2n+1}
\end{array}
\right.
$$

Values of $U_i$ and $V_i$ in previous tables can be computed easily
by the following PARI/gp programs:

$U_{2j+1}$ : \texttt{U0=1;U1=6; for(i=1,N, U0=5*U1-U0; U1=5*U0-U1;
print(4*i+1," ",U0); print(4*i+1," ",U1))}

\end{document}